\documentclass[a4paper, 10pt]{article}

\usepackage{geometry}
\usepackage[francais, english]{babel}
\usepackage{latexsym,amsfonts, amssymb}
\usepackage{amsthm}
\usepackage[latin1]{inputenc}
\usepackage{amsmath}
\usepackage{amssymb}
\usepackage{amsfonts}
\usepackage{dsfont}
\usepackage{color}
\usepackage{multicol}
\usepackage{hyperref}
\usepackage{amsthm}
\usepackage{latexsym}
\usepackage{fancyhdr}
\usepackage{epsfig}
\usepackage{mathrsfs}
\usepackage{multirow}

\makeatletter
\renewcommand\theequation{\thesection.\arabic{equation}}
\@addtoreset{equation}{section} \makeatother

\newfont{\bbf}{cmbx12 scaled 1435}

\newtheorem{thm}{Theorem}[section]

\newtheorem{lem}{Lemma}[section]

\renewcommand{\theequation}{\thesection.\arabic{equation}}

\begin{document}

\newcommand{\I}[1]{\mathds{1}_{{#1}}}

\def\Sum{ \displaystyle \sum }
\def\Frac{\displaystyle \frac}

\def\Cit{\mathbb{C}}
\def\esp{\mathbb{E}}
\def\Var{\hbox{\rm Var}}
\def\Cov{\hbox{\rm Cov}}
\def\Supp{\hbox{\rm Supp}}
\def\Card{\hbox{\rm Card}}
\def\Det{\hbox{\rm Det}}
\def\Tr{\hbox{\rm Tr}}
\def\Dim{\hbox{\rm dim}}
\def\Rank{\hbox{\rm dim}}
\def\Id{\hbox{\rm Id}}
\def\Ker{\hbox{\rm Ker}}
\def\ind{\mathbb{I}}
\def\Nit{\mathbb{N}}
\def\Rit{\mathbb{R}}
\def\Zit{\mathbb{Z}}
\def\prob{\mathbb{P}}
\def\where{\rm where}
\def\with{\rm with}
\def\I{\rm I}
\def\J{\rm J}
\def\Ip{\rm I^{\prime}}
\def\Jp{\rm J^{\prime}}
\def\as{\rm a.s}

\title{Estimation of the Error Density  in a \\Semiparametric Transformation Model}

\author{Rawane \textsc{Samb}
\\
Universit\'e catholique de Louvain
\thanks{R. Samb acknowledges financial support from IAP research network P6/03
of the Belgian Government (Belgian Science Policy).} \\
\and
C\'edric \textsc{Heuchenne} \\
University of Li\`ege and Universit\'e catholique de Louvain
\thanks{C. Heuchenne  acknowledges financial support from IAP research network P6/03 of the Belgian Government (Belgian Science Policy), and from the contract `Projet d'Actions de Recherche Concert\'ees' (ARC) 11/16-039 of the `Communaut\'e fran\c{c}aise de Belgique', granted by the `Acad\'emie universitaire Louvain'.} \\
\and
Ingrid \textsc{Van Keilegom} \\
Universit\'e catholique de Louvain
\thanks{I. Van Keilegom acknowledges financial support from IAP research network P6/03 of the Belgian Government (Belgian Science Policy), from the European Research Council under the European Community's Seventh Framework Programme (FP7/2007-2013) / ERC Grant agreement No. 203650, and from the contract `Projet d'Actions de Recherche Concert\'ees' (ARC) 11/16-039 of the `Communaut\'e fran\c{c}aise de Belgique', granted by the `Acad\'emie universitaire Louvain'.} \\}

\date{\today}

\maketitle

\begin{abstract}
Consider the semiparametric transformation model
$\Lambda_{\theta_o}(Y)=m(X)+\varepsilon$, where $\theta_o$ is an
unknown finite dimensional parameter, the functions
$\Lambda_{\theta_o}$ and $m$ are smooth, $\varepsilon$ is
independent of $X$, and $\esp(\varepsilon)=0$. We propose a kernel-type estimator of the density of the error
$\varepsilon$, and prove its asymptotic normality.
The estimated errors, which lie at the basis of this estimator, are obtained
from a profile likelihood estimator of
$\theta_o$ and a nonparametric kernel estimator of $m$. The
practical performance of the proposed density estimator is evaluated in a simulation study. \\
\end{abstract}


\newpage
\normalsize

\renewcommand{\thefootnote}{\arabic{footnote}}
\setcounter{footnote}{1} \setlength{\baselineskip}{.26in}

\section{Introduction}

Let $(X_{1},Y_{1}),\ldots,(X_{n},Y_{n})$ be independent replicates
of the random vector $(X, Y)$, where $Y$ is a univariate
dependent variable and $X$ is a one-dimensional covariate. We
assume that $Y$ and $X$ are related via the semiparametric
transformation model
\begin{equation}
\Lambda_{\theta_o}(Y) = m(X) + \varepsilon,  \label{model}
\end{equation}
where $\varepsilon$ is independent of $X$ and has mean zero. We
assume that $\{\Lambda_{\theta}: \theta\in\Theta\}$ (with $\Theta \subset \Rit^p$ compact) is a
parametric family of strictly increasing functions defined on an
unbounded subset $\mathcal{D}$ in $\Rit$, while $m$ is the unknown
regression function, belonging to an infinite dimensional parameter set $\mathcal{M}$. We assume
that $\mathcal{M}$ is a space of functions endowed with the
norm $\|\cdot\|_{\mathcal{M}}=\|\cdot\|_{\infty}$. We
denote $\theta_o\in\Theta$ and $m\in\mathcal{M}$ for the true
unknown finite and infinite dimensional parameters. Define the
regression function
$$
 m_{\theta}(x) = \esp[\Lambda_{\theta}(Y)|X=x],
$$
for each $\theta\in\Theta$, and let
$\varepsilon_{\theta}=\varepsilon(\theta)=\Lambda_{\theta}(Y)-m_{\theta}(X)$.

In this paper, we are interested in the estimation of the
probability density function (p.d.f.) $f_{\varepsilon}$  of the
residual term $\varepsilon=\Lambda_{\theta_o}(Y)-m(X)$. To this end, we
first obtain the estimators $\widehat{\theta}$ and
$\widehat{m}_{\theta}$ of the parameter $\theta_o$ and the function
$m_{\theta}$, and second, form the semiparametric regression
residuals $\widehat{\varepsilon}_i(\widehat{\theta})
=\Lambda_{\widehat{\theta}}(Y_i)-\widehat{m}_{\widehat{\theta}}(X_i)$.
To estimate $\theta_o$ we use a profile likelihood (PL) approach,
developed in Linton, Sperlich and
Van Keilegom (2008), whereas $\widehat{m}_{\theta}$ is estimated by
means of a Nadaraya-Watson-type estimator (Nadaraya, 1964, Watson, 1964).
To our knowledge, the
estimation of the density of $\varepsilon$ in model
(\ref{model}) has not yet been investigated in the statistical
literature. This estimation may be very useful in various
regression problems. Indeed, taking transformations of the data may
induce normality and error variance homogeneity in the transformed model. So the
estimation of the error density in the transformed model may
be used for testing these hypotheses.

Taking transformations of
the data has been an important part of statistical practice for
many years. A major contribution to this methodology was made by
Box and Cox (1964), who proposed a parametric power family of
transformations that includes the logarithm and the identity. They
suggested that the power transformation, when applied to the
dependent variable in a linear regression model, might induce normality
and homoscedasticity. Lots of
effort  has been devoted to the investigation of the Box-Cox
transformation since its introduction. See, for example,
 Amemiya (1985), Horowitz (1998), Chen,
Lockhart and Stephens (2002), Shin (2008), and Fitzenberger, Wilke and
Zhang (2010).   Other dependent variable transformations have
been suggested, for example, the Zellner and Revankar (1969)
transform and the Bickel and Doksum (1981) transform. The
transformation methodology has been quite successful and a large
literature exists on this topic for parametric models. See
Carroll and Ruppert (1988) and Sakia (1992)  and references
therein.

The estimation of (functionals of) the error distribution and density under
simplified versions of model (\ref{model}) has received considerable
attention in the statistical literature in recent years.  Consider e.g.\ model
(\ref{model}) but with $\Lambda_{\theta_o} \equiv id$, i.e.\ the response is
not transformed. Under this model, Escanciano and Jacho-Chavez (2010)
considered the estimation of the (marginal) density of the response $Y$ via
the estimation of the error density.  Akritas and Van Keilegom (2001)
estimated the cumulative distribution function of the regression
error in a heteroscedastic model with univariate covariates. The
estimator they proposed is based on nonparametrically estimated regression
residuals.  The weak convergence of their estimator  was proved. The results
obtained by Akritas and Van Keilegom (2001) have been generalized by
Neumeyer and Van Keilegom (2010) to the case
of multivariate covariates. M\"{u}ller, Schick and Wefelmeyer
(2004) investigated linear functionals of the error distribution in
nonparametric regression. Cheng (2005) established the asymptotic
normality of an estimator of the error density based on estimated
residuals. The estimator he proposed is constructed by splitting
the sample into two parts: the first part is used for the
estimation of the
 residuals, while the second part of the sample is used for the
construction of the error density estimator.  Efromovich (2005)
proposed an adaptive estimator of the error density, based on a
density estimator proposed by Pinsker (1980). Finally, Samb (2010)
also considered the estimation of the error density, but his
approach is more closely related to the one in Akritas and Van
Keilegom (2001).

In order to achieve the objective of this paper, namely the
estimation of the error density under model (\ref{model}), we
first need to estimate the transformation parameter $\theta_o$.
To this end, we make use of the results in Linton, Sperlich and
Van Keilegom (2008). In the latter paper, the authors first
discuss the nonparametric identification of model (\ref{model}),
and second, estimate the transformation parameter $\theta_o$ under
the considered model. For the estimation of this parameter, they
propose two approaches. The first approach uses a semiparametric
profile likelihood (PL) estimator, while the second is based on a
`mean squared distance from independence-estimator (MD) using the
estimated distributions of $X$, $\varepsilon$ and
$(X,\varepsilon)$.  Linton, Sperlich and Van Keilegom (2008)
derived the asymptotic distributions of their estimators under
certain regularity conditions, and proved that both estimators of
$\theta_o$ are root-$n$ consistent. The authors also showed that,
in practice, the performance of the PL method is better than that
of the MD approach. For this reason, the PL method will be
considered in this paper for the estimation of $\theta_o$.

The rest of the paper is organized as follows. Section 2 presents
our estimator of the error density and groups some notations and technical assumptions.
Section 3 describes the asymptotic results of the
paper.  A simulation study is given in Section 4, while
Section 5 is devoted to some general conclusions.
Finally, the proofs of the asymptotic results are collected in Section 6.

\section{Definitions and assumptions}

\subsection{Construction of the estimators}

The approach proposed here for the estimation of $f_{\varepsilon}$
is based on a two-steps procedure. In a first step, we estimate
the finite dimensional parameter $\theta_o$. This parameter is
estimated by the profile likelihood (PL) method, developed in
Linton, Sperlich and Van Keilegom (2008). The
basic idea of this method is  to replace all unknown expressions
in the likelihood function by their nonparametric kernel
estimates. Under model (\ref{model}), we have
$$
\prob\left(Y\leq y|X\right)
= \prob\left(\Lambda_{\theta_o}(Y) \leq \Lambda_{\theta_o}(y)|X \right)
= \prob \left(\varepsilon_{\theta_o} \leq \Lambda_{\theta_o}(y) - m_{\theta_o}(X)|X \right)
= F_{\varepsilon} \left( \Lambda_{\theta_o}(y)-m_{\theta_o}(X) \right).
$$
Here, $F_{\varepsilon}(t)=\prob(\varepsilon\leq t)$, and so
$$
f_{Y|X}(y|x)
= f_{\varepsilon} \left( \Lambda_{\theta_o}(y)-m_{\theta_o}(x) \right)
\Lambda_{\theta_o}^{\prime}(y),
$$
where $f_{\varepsilon}$ and $f_{Y|X}$ are the densities
 of $\varepsilon$, and of $Y$ given $X$, respectively. Then, the
$\log$ likelihood function is
$$
\sum_{i=1}^n
\left\lbrace
\log f_{\varepsilon_{\theta}}
(\Lambda_{\theta}(Y_i)-m_{\theta}(X_i))
 +
\log \Lambda_{\theta}^{\prime}(Y_i)
\right\rbrace,
\quad\theta\in\Theta,
$$
where $f_{\varepsilon_{\theta}}$ is the density function of
$\varepsilon_{\theta}$. Now, let
\begin{equation}
\widehat{m}_{\theta}(x)
=
\frac{
\sum_{j=1}^n
\Lambda_{\theta}(Y_j)
 K_1\left(\frac{X_j-x}{h}\right)}
  {\sum_{j=1}^n
  K_1\left(\frac{X_j-x}{h}\right)}
   \label{NW}
\end{equation}
be the Nadaraya-Watson estimator of $m_{\theta}(x)$, and let
\begin{equation}
 \widehat{f}_{\varepsilon_{\theta}}(t)
 =
 \frac{1}{ng}
\sum_{i=1}^n
K_2 \left(
\frac{\widehat{\varepsilon}_i(\theta)-t}{g}
\right).
\label{DE}
\end{equation}
where
 $\widehat{\varepsilon}_i(\theta)
=\Lambda_{\theta}(Y_i)-\widehat{m}_{\theta}(X_i)$. Here, $K_1$ and
$K_2$ are kernel functions and $h$ and $g$ are bandwidth
sequences. Then, the PL estimator of $\theta_o$ is defined by
\begin{equation}
\widehat{\theta}
=
\arg\max_{\theta\in\Theta}
\sum_{i=1}^n
\left[
\log \widehat{f}_{\varepsilon_{\theta}}
(\Lambda_{\theta}(Y_i)-\widehat{m}_{\theta}(X_i))
+
 \log\Lambda_{\theta}^{\prime}(Y_i)
 \right].
 \label{PLE}
\end{equation}
\textcolor{black}{Recall that $\widehat{m}_{\theta}(X_i)$ converges to $m_\theta(X_i)$ at a slower rate for
those $X_i$ which are close to the boundary of the support $\mathcal{X}$ of the covariate $X$.}
 That is why we assume implicitly that the proposed estimator
 (\ref{PLE}) of $\theta_o$ trims the observations $X_i$
 outside \textcolor{black}{a subset $\mathcal{X}_0$  of $\mathcal{X}$}. Note that we keep the root-$n$ consistency of
$\widehat{\theta}$ proved in Linton, Sperlich and Van Keilegom
(2008) by trimming the covariates outside $\mathcal{X}_0$. But in
this case, the resulting asymptotic variance is different to the
one obtained in the latter paper.

\noindent In a second step, we use the above estimator
$\widehat{\theta}$ to build the estimated residuals
$\widehat{\varepsilon}_i(\widehat{\theta})
=\Lambda_{\widehat{\theta}}(Y_i)-\widehat{m}_{\widehat{\theta}}(X_i)$.
Then, our proposed estimator
$\widehat{f}_{\widehat{\varepsilon}}(t)$ of $f_{\varepsilon}(t)$
is defined by
\begin{equation}
\widehat{f}_{\widehat{\varepsilon}}(t)
=
\frac{1}{nb} \sum_{i=1}^n
K_3\left(
\frac{\widehat{\varepsilon}_i(\widehat{\theta})-t}{b}
\right),
 \label{pdf}
\end{equation}
where $K_3$ is a kernel function and $b$ is a bandwidth sequence,
not necessarily the same as the kernel $K_2$ and the bandwidth $g$ used in (\ref{DE}).
Observe that this estimator  is a feasible estimator in the sense
that it does not depend on any unknown
 quantity, as is desirable in practice. This contrasts with the unfeasible
 ideal kernel estimator
\begin{equation}
\widetilde{f}_{\varepsilon}(t)
=
\frac{1}{nb}
\sum_{i=1}^n
K_3\left(\frac{\varepsilon_i-t}{b}\right),
 \label{ppdf}
\end{equation}
 which depends in particular on the unknown regression errors
$\varepsilon_i=\varepsilon_i(\theta_o)=\Lambda_{\theta_o}(Y_i)-m(X_i)$.
 It is however intuitively clear that $\widehat{f}_{\widehat{\varepsilon}}(t)$ and
$\widetilde{f}_{\varepsilon}(t)$ will be very close for $n$ large enough,
as will be illustrated by the results given in Section 3.

\subsection{Notations}

When there is no ambiguity, we use $\varepsilon$ and $m$ to
indicate $\varepsilon_{\theta_o}$ and $m_{\theta_o}$. Moreover,
$\mathcal{N}(\theta_o)$ represents a neighborhood of $\theta_o$.
For the kernel $K_j$ $(j=1,2,3)$,  let $\mu(K_j) = \int v^2 K_j(v)dv$ and let
 $K_j^{(p)}$ be the $p{\rm th}$ derivative of $K_j$. For any function
$\varphi_{\theta}(y)$, denote
$\dot{\varphi}_{\theta}(y)=\partial\varphi_{\theta}(y)/\partial\theta
= (\partial\varphi_{\theta}(y)/\partial\theta_1,\ldots,\partial\varphi_{\theta}(y)/\partial\theta_p)^t$
and $\varphi_{\theta}^{\prime}(y)
=\partial\varphi_{\theta}(y)/\partial y$.
Also, let $\|A\| = (A^t A)^{1/2}$ be the Euclidean norm of any vector $A$.

 For any functions $\widetilde{m}$, $r$, $f$, $\varphi$ and $q$, and
 any $\theta\in\Theta$,
  let
 $s=(\widetilde{m},r,f,\varphi,q)$,
$s_{\theta}=(m_{\theta},\dot{m}_{\theta},f_{\varepsilon_{\theta}},
f_{\varepsilon_{\theta}}^{\prime},
\dot{f}_{\varepsilon_{\theta}})$,
$\varepsilon_i(\theta,\widetilde{m})
=\Lambda_{\theta}(Y_i)-\widetilde{m}(X_i)$, and define
$$
G_n(\theta,s)
=
n^{-1}
\sum_{i=1}^n
\left\lbrace
\frac{1}{f\{\varepsilon_i(\theta,\widetilde{m})\}}
\left[
\varphi\{\varepsilon_i(\theta,\widetilde{m})\}
\{
\dot{\Lambda}_{\theta}(Y_i) - r(X_i) \}
+
q\{\varepsilon_i(\theta,\widetilde{m})\}
\right]
+
\frac{\dot{\Lambda}_{\theta}^{\prime}(Y_i)}
{\Lambda_{\theta}^{\prime}(Y_i)}
\right\rbrace,
$$
$G(\theta,s)=\esp[G_n(\theta,s)]$ and
$\mathcal{G}(\theta_o,s_{\theta_o}) =
\left.\frac{\partial}{\partial \theta}G(\theta,s_\theta)
\right\downarrow_{\theta=\theta_o}$.

\subsection{Technical assumptions}

The assumptions we need for the asymptotic results are
listed below for convenient reference.

\vskip 0.1cm \noindent {$\bf(A1)$} The function $K_j$ $(j=1,2,3)$
is symmetric, has compact support, $\int\! v^k K_j(v) dv=0$ for
$k=1,\ldots,q_{j}-1$ and $\int\!v^{q_j} K_j(v) dv\neq 0$ for some
$q_j \ge 4$, $K_j$ is twice continuously differentiable, and
$\int\! K_3^{(1)}(v) dv=0$.
\medskip
\\
{$\bf(A2)$} The bandwidth sequences $h$, $g$ and $b$ satisfy
$nh^{2q_1}=o(1)$, $ng^{2q_2}=o(1)$
(where $q_1$ and $q_2$ are
defined in (A1)), $(nb^5)^{-1} = O(1)$,
$nb^3h^2(\log h^{-1})^{-2}\rightarrow\infty$
and $ng^6(\log g^{-1})^{-2}\rightarrow\infty$.
\medskip
\\
{$\bf(A3)$} (i) The support $\mathcal{X}$ of the covariate $X$ is a
compact subset of $\Rit$, and $\mathcal{X}_0$
 is a subset with non empty interior, whose closure is in the
  interior of $\mathcal{X}$.
\\
(ii) The density $f_X$ is bounded away from zero and infinity on
${\cal X}$, and has continuous second order partial derivatives on
${\cal X}$.
\medskip
\\
{$\bf(A4)$} The function $m_\theta(x)$ is twice continuously
differentiable with respect to $\theta$ on ${\cal X} \times {\cal
N}(\theta_0)$, and the functions $m_{\theta}(x)$ and
$\dot{m}_{\theta}(x)$ are $q_1$ times continuously differentiable
with respect to $x$ on $\mathcal{X} \times {\cal N}(\theta_0)$.
All these derivatives are bounded, uniformly in
$(x,\theta)\in\mathcal{X}\times\mathcal{N}(\theta_o)$.
\medskip
\\
{$\bf(A5)$}
The error $\varepsilon=\Lambda_{\theta_o}(Y)-m(X)$
has finite fourth moment and is independent of $X$.
\medskip
\\
{$\bf(A6)$}
The distribution $F_{\varepsilon_\theta}(t)$ is $q_3+1$ (respectively three) times
continuously differentiable with respect to $t$ (respectively $\theta$), and
$$
\sup_{\theta,t} \left\| \frac{\partial^{k+\ell}}
{\partial t^k\partial\theta_1^{\ell_1}\ldots\partial\theta_p^{\ell_p}}
F_{\varepsilon_\theta}(t) \right\| <\infty
$$
for all $k$ and $\ell$ such that $0\leq k+\ell\leq
2$, where $\ell=\ell_1+\ldots+\ell_p$ and
$\theta=(\theta_1,\ldots,\theta_p)^t$.
\medskip
\\
{$\bf(A7)$}
The transformation $\Lambda_{\theta}(y)$ is three times
continuously differentiable with respect to both $\theta$ and $y$, and there
exists a $\alpha>0$ such that
$$
\esp\left[\sup_{\theta^\prime : \|\theta^{\prime}-\theta\|\leq\alpha}
\left\|\frac{\partial^{k+\ell}}
{\partial y^k\partial\theta_1^{\ell_1}\ldots\partial\theta_p^{\ell_p}}
\Lambda_{\theta^{\prime}}(Y) \right\| \right] <\infty
$$
\vskip 0.15cm\noindent for all $\theta\in\Theta$, and for all $k$
and $\ell$ such that $0\leq k+\ell\leq 3$, where
$\ell=\ell_1+\ldots+\ell_p$ and
$\theta=(\theta_1,\ldots,\theta_p)^t$. Moreover, $\sup_{x \in
{\cal X}} \esp[\dot\Lambda_{\theta_o}^4(Y)|X=x] < \infty$.
\medskip
\\
{$\bf(A8)$}
For all $\eta>0$, there exists $\epsilon(\eta)>0$ such that
$$
\inf_{\|\theta-\theta_o\|>\eta} \|G(\theta,s_{\theta})\| \geq \epsilon(\eta)>0.
$$
Moreover, the matrix $\mathcal{G}(\theta_o,s_{\theta_o})$ is non-singular.
\medskip
\\
{$\bf(A9)$}
(i) $\esp(\Lambda_{\theta_o}(Y))=1$, $\Lambda_{\theta_o}(0)=0$ \textcolor{black}{and the set $\{x \in {\cal X}_0 : m'(x) \neq 0\}$ has nonempty interior.}
\medskip
\\
(ii) Assume that $\phi(x,t)
=\dot{\Lambda}_{\theta_o}
(\Lambda_{\theta_o}^{-1}(m(x)+t)) f_{\varepsilon}(t)$ is continuously differentiable
 with respect to $t$ for all $x$ and that
\begin{equation}
\sup_{s:|t-s|\leq\delta}
\esp
\left|
\frac{\partial \phi}
{\partial s} (X,s)
\right|
<\infty.
 \label{cond}
\end{equation}
for all $t \in \Rit$ and for some $\delta>0$.

\vskip 0.3cm
Assumptions (A1), part of (A2), (A3)(ii), (A4) and (A6), (A7) and (A8) are
 used by Linton, Sperlich and Van Keilegom (2008) to show that
 the PL estimator $\widehat{\theta}$ of $\theta_o$ is root
 $n$-consistent. The differentiability of $K_j$ up to second order imposed in
assumption (A1) is used to expand the
 two-steps kernel estimator
$\widehat{f}_{\widehat{\varepsilon}}(t)$ in (\ref{pdf}) around the
unfeasible one $\widetilde{f}_{\varepsilon}(t)$. Assumptions
(A3)(ii) and (A4) impose that all the functions to be estimated
have  bounded derivatives. The last assumption in (A2) is useful
for obtaining the uniform convergence of the Nadaraya-Watson
estimator of $m_{\theta_o}$ in (\ref{NW}) (see for instance
Einmahl and Mason, 2005). This  assumption  is also needed in the
study of  the difference between the feasible estimator
$\widehat{f}_{\widehat{\varepsilon}}(t)$ and the unfeasible
estimator $\widetilde{f}_{\varepsilon}(t)$.  Finally, (A9)(i) is
needed for identifying the model (\textcolor{black}{see Vanhems and
Van Keilegom (2011)}).

\section{Asymptotic results}

In this section we are interested in the asymptotic behavior of
the estimator $\widehat{f}_{\widehat{\varepsilon}}(t)$. To this
end, we first investigate its asymptotic representation,  which will
be needed to show its asymptotic normality.

\begin{thm}
Assume {\rm (A1)-(A9)}.  Then,
$$
\widehat{f}_{\widehat{\varepsilon}}(t)
-
f_{\varepsilon}(t)
=
\frac{1}{nb}
\sum_{i=1}^n
K_{3}\left(\frac{\varepsilon_i-t}{b}\right)
-
 f_{\varepsilon}(t)
+
R_n(t),
$$
where $R_n(t)=o_{\prob}\left((n b)^{-1/2}\right)$ for all $t\in\Rit$.
\label{Asymp-Dev}
\end{thm}

This result is important, since it shows that,
provided the bias term is negligible, the estimation of $\theta_o$ and
$m(\cdot)$ has asymptotically no effect on the behavior of the
 estimator $\widehat{f}_{\widehat{\varepsilon}}(t)$.
Therefore, this estimator is asymptotically equivalent to the unfeasible
 estimator $\widetilde f_{\varepsilon}(t)$, based on
the unknown true errors $\varepsilon_1,\ldots,\varepsilon_n$.

Our next result gives the asymptotic normality of the estimator
$\widehat{f}_{\widehat{\varepsilon}}(t)$.

\begin{thm}
Assume {\rm (A1)-(A9)}.  In addition, assume that
$nb^{2q_3+1}=O(1)$.
Then,
$$
\sqrt{nb}
\left(
\widehat{f}_{\widehat{\varepsilon}}(t)
 -
\overline{f}_{\varepsilon}(t)
\right)
\stackrel{d}{\rightarrow}
N\left(0,f_\varepsilon(t)
\int K_3^2(v) dv \right),
$$
where
$$
\overline{f}_{\varepsilon}(t)
 =
 f_{\varepsilon}(t)
 +
\frac{b^{q_3}}{q_3!}
f_{\varepsilon}^{(q_3)}(t)
\int
v^{q_3}
K_3(v)
dv.
$$
\label{Normality}
\end{thm}

\noindent The proofs of Theorems \ref{Asymp-Dev} and \ref{Normality}
are given in Section \ref{proofs}.

\section{Simulations}

In this section, we investigate the performance of our method for
different models and different sample sizes.
Consider
\begin{equation}
\Lambda_{\theta_o}(Y)
=
b_0+b_1X^2+b_2\sin(\pi X)
+
\sigma_{e} \varepsilon,
\label{model2}
\end{equation}
where $\Lambda_{\theta}$ is the Manly (1976) transformation
$$
\Lambda_{\theta}(y)
=
\begin{cases}
\frac{e^{\theta y}-1}{\theta}, & \theta\neq 0, \\
y, & \theta=0,
\end{cases}
$$
$\theta_o\in[-0.5,1.5]$,  $X$ is uniformly distributed
on the interval $[-0.5,0.5]$,
and $\varepsilon$ is independent of $X$ and has a standard
normal distribution but restricted
to the interval $[-3,3]$. We study three different model settings.
 For each of them, $b_0=3\sigma_{e}+b_2$.
The other parameters are chosen as follows: \\[.3cm]
\begin{tabular}{llll}
\hspace*{1cm} {\rm Model 1:} & $b_1=5$, & $b_2=2$, & $\sigma_e=1.5$; \\
\hspace*{1cm} {\rm Model 2:} & $b_1=3.5$, & $b_2=1.5$, & $\sigma_e=1$; \\
\hspace*{1cm} {\rm Model 3:} & $b_1=2.5$, & $b_2=1$, & $\sigma_e=0.5$. \\[.3cm]
\end{tabular}

The parameters and the error distribution have been chosen in such
a way that the transformation $\Lambda_{\theta_o}(Y)$ is positive,
to avoid problems when generating the variable $Y$. Our
simulations are done for $\theta_o = 0$, $0.5$ or $1$. The
estimator of $\theta_o$ is chosen from a grid on the interval
$[-0.5,1.5]$ with step size $0.0625$. We used the kernel $K(x) =
\frac{15}{16}\left(1-x^2\right)^2 \mathds{1}\left(|x|\leq
1\right)$ for both the regression function and
 the density estimators. The results are based on $100$ random
samples of size $n=50$ or $n=100$, and we worked with the bandwidths
$h=0.3\times n^{-1/5}$ and $b=g=r_n$, where $r_n=1.06\times{\rm
std}(\widehat{\varepsilon})\times n^{-1/5}$, which is Silverman's
(1986) rule of thumb bandwidth for univariate density estimation.
Here ${\rm std}(\widehat{\varepsilon})$ is the average
of the standard deviations of $\widehat\varepsilon$ over the 100 samples.

\textcolor{black}{Table \ref{tab1} shows the values of the
mean, standard deviation and mean squared error of
$\widehat{\theta}$
for the considered models, sample sizes and values of
$\theta_o$.  We observe that the results for the different models
are quite similar, and as expected, the results are better for $n=100$ than for $n=50$. }

\begin{figure}[H]
 \centering
\includegraphics[width=14cm, height=15cm]{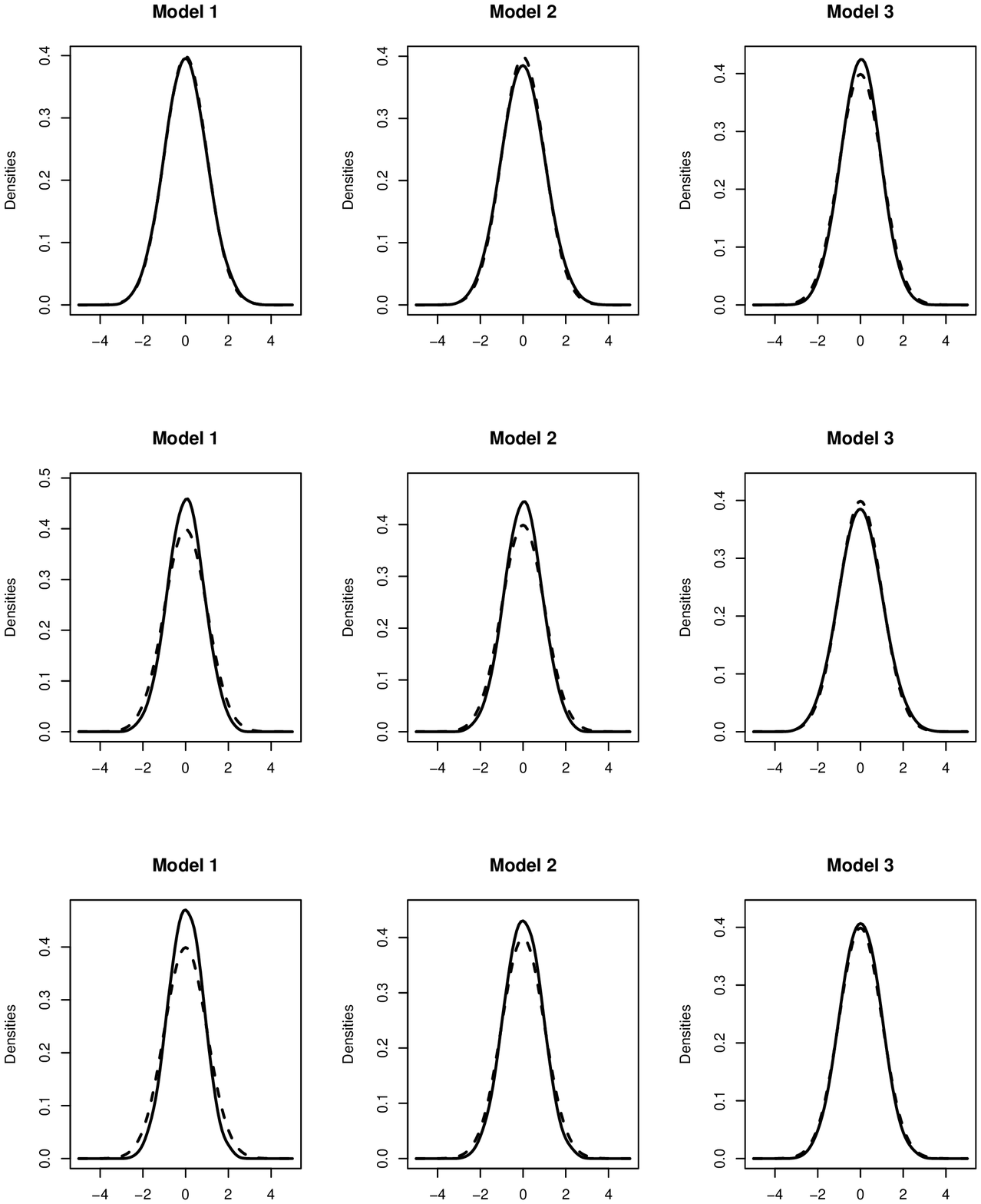}
\caption{Curves of the pointwise average of $\widehat{f}_{\widetilde\varepsilon}$
over $100$ random samples of size $n=100$ (solid curve) and
of the standard normal density (dashed curve) for $\theta_o=0$ (first row),
$\theta_o=0.5$ (second row) and $\theta_o=1$ (third row).}
\label{Fig1}
\end{figure}

Table \ref{tab2} shows the mean squared error (MSE) of the
estimator $\widehat{f}_{\widetilde\varepsilon}(t)$ of the
standardized (pseudo-estimated) error $\widetilde \varepsilon =
(\Lambda_{\widehat\theta}(Y)-\widehat{m}_{\widehat\theta}(X))/\sigma_e$,
for sample sizes $n=50$ and $n=100$ and for $t=-1$, $0$ and $1$.
Results for $\widehat f_{\widehat \varepsilon}(t)$ have also been
obtained, but are not reported here.  Indeed, Figure \ref{Fig1},
displaying $\widehat f_{\widetilde \varepsilon}(t)$, shows that,
even though residuals are standardized for each simulation (with
known $\sigma_e$), better behavior is observed for models with
smaller $\sigma_e$. Moreover, we observe that for $\theta_o=0$
there is very little difference between the curve of
$\widehat{f}_{\widetilde\varepsilon}$ and the one of the standard
normal density. On the other hand for $\theta_o=0.5$ and
$\theta_o=1$, we notice an important difference between the two
curves under Model 1 and 2, but the difference is less important
under Model 3.

\begin{table}[h]
\begin{center}
\renewcommand{\arraystretch}{1.1}
\begin{tabular}[b]{|c|c|c|c|c|c|c|c|c|c|c|}
\hline $n$& $\theta_o$ &
\multicolumn{3}{|c|}{mean($\widehat{\theta}$)} &
\multicolumn{3}{|c|}{std($\widehat{\theta}$)}  &
\multicolumn{3}{|c|}{MSE($\widehat{\theta}$)}
 \\
\cline{3-11}
&  & Model 1 & Model 2 & Model 3  & Model 1 & Model 2 & Model 3 & Model 1 & Model 2 & Model 3
\\
\hline 50
& 0 & 0.0063 & 0.0065 & 0.0071  & 0.0116 & 0.0161 & 0.0239 & 0.0064 & 0.0124 & 0.0277\\
& 0.5 & 0.3787 & 0.3754 & 0.3907 & 0.0417 & 0.0438 & 0.0486 & 0.0783 & 0.0867 & 0.1140\\
& 1 & 0.8197 & 0.8449 & 0.8658 & 0.0792 & 0.0796 & 0.0798 & 0.3506 & 0.3492 & 0.3409 \\
\hline
100 & 0  & 0.0055 & 0.0148 & 0.0170 & 0.0057 & 0.0078 & 0.0116 & 0.0032 & 0.0059 & 0.0132 \\
& 0.5 & 0.4596 & 0.4621 & 0.4728 & 0.0246 & 0.254 & 0.0270 & 0.0634 & 0.0676 & 0.0752 \\
& 1    & 0.9196 & 0.9545 & 0.9749 & 0.0401 & 0.0437 & 0.0438 & 0.2092 & 0.1999 & 0.1637 \\
\hline
\end{tabular}
\end{center}
\caption{Approximation of the mean,  the
standard deviation and the mean squared error of $\widehat{\theta}$
for the three regression
models. All numbers are calculated based on $100$ random
samples.}
\label{tab1}
\end{table}

\begin{table}[h]
\begin{center}
\renewcommand{\arraystretch}{1.3}
\begin{tabular}[b]{|c|c|c|c|c|c|}
\hline $n$& $\theta_o$ & $t$ & \multicolumn{3}{|c|}{Mean Squared
Error of $\widehat{f}_{\widetilde\varepsilon}(t)$}
\\
\cline{4-6} &  &  &Model 1 & Model 2 & Model 3 \\
\hline
50
& 0 &-1 & 0.0026 & 0.0025 & 0.0019 \\
&    & 0 & 0.0040 & 0.0037 & 0.0028 \\
&    & 1 & 0.0026 & 0.0023 & 0.0017 \\
\cline{2-6}
& 0.5 & -1 & 0.0063 & 0.0048 & 0.0025 \\
&       & 0   & 0.0527 & 0.0372 & 0.0147 \\
&       & 1   & 0.0062 & 0.0046 & 0.0020 \\
\cline{2-6}
& 1    & -1 & 0.0078 & 0.0048 & 0.0024 \\
&       & 0  & 0.0564 & 0.0314 & 0.0133 \\
&       & 1 & 0.0049 & 0.0030 & 0.0019 \\
\hline
100
& 0 & -1 & 0.0012 & 0.0011 & 0.0008 \\
&    & 0  & 0.0039 & 0.0035 & 0.0026 \\
&    & 1  & 0.0017 & 0.0015 & 0.0012 \\
\cline{2-6}
&0.5& -1& 0.0015 & 0.0014 & 0.0011 \\
&     & 0  & 0.0075 & 0.0057 & 0.0031 \\
&     & 1  & 0.0021 & 0.0018 & 0.0012 \\
\cline{2-6}
& 1  & -1& 0.0024 & 0.0018 & 0.0012 \\
&     & 0  & 0.0110 & 0.0052 & 0.0270 \\
&     & 1  & 0.0019 & 0.0016 & 0.0012 \\
\hline
 \end{tabular}
\end{center}
 \caption{Mean squared error of
 $\widehat{f}_{\widetilde\varepsilon}(t)$ for three regression models.
All numbers are calculated based on $100$ random samples.}
 \label{tab2}
\end{table}

\section{Conclusions}

In this paper we have studied the estimation of the density of the error
in a semiparametric transformation model.  The regression function in this model
 is unspecified (except for some smoothness assumptions), whereas the transformation
 (of the dependent variable in the model) is supposed to belong to a parametric family
 of monotone transformations.  The proposed estimator is a kernel-type estimator,
 and we have shown its asymptotic normality.  The finite sample performance of the
  estimator is illustrated by means of a simulation study.

It would be interesting to explore various possible applications of the results in this paper.
For example, one could use the results on the estimation of the error density to test hypotheses
concerning e.g.\ the normality of the errors, the homoscedasticity of the error variance,
or the linearity of the regression function, all of which are important features in the context
of transformation models.

\section{Proofs} \label{proofs}

\noindent
{\bf Proof of Theorem \ref{Asymp-Dev}.} Write
$$
\widehat{f}_{\widehat{\varepsilon}}(t)
-
 f_{\varepsilon}(t)
  =
[\widehat{f}_{\varepsilon}(t) - f_{\varepsilon}(t)]
+
[\widehat{f}_{\widehat{\varepsilon}}(t)
 -
\widehat{f}_{\varepsilon}(t)],
$$
where
$$
\widehat{f}_{\varepsilon}(t)
=
 \frac{1}{nb}
 \sum_{i=1}^n
 K_3\left(\frac{\widehat{\varepsilon}_i-t}{b}\right)
$$
and $\widehat{\varepsilon}_i
=\Lambda_{\theta_o}(Y_i)-\widehat{m}_{\theta_o}(X_i)$,
$i=1,\ldots,n$. In a completely similar way as was done for Lemma
A.1 in Linton, Sperlich and Van Keilegom (2008), it can be shown
that
\begin{eqnarray}
\label{rem}
\widehat{f}_{\varepsilon}(t)
-
 f_{\varepsilon}(t)
=
\frac{1}{nb}
 \sum_{i=1}^n
 K_3\left(\frac{\varepsilon_i-t}{b}\right)
 -
f_{\varepsilon}(t)
+
o_{\prob}((n b)^{-1/2})
\end{eqnarray}
 for all $t\in\Rit$. Note that the remainder term in Lemma A.1 in
the above paper equals a sum of i.i.d.\ terms of mean zero, plus a
$o_{\prob}(n^{-1/2})$ term.  Hence, the remainder term in that
paper is $O_{\prob}(n^{-1/2})$, whereas we write
$o_{\prob}((nb)^{-1/2})$ in (\ref{rem}).
 Therefore, the result of the theorem follows if we prove that $
\widehat{f}_{\widehat{\varepsilon}}(t)
-\widehat{f}_{\varepsilon}(t) = o_{\prob}((nb)^{-1/2})$. To this
end, write
\begin{eqnarray}
\nonumber
\lefteqn{
\widehat{f}_{\widehat{\varepsilon}}(t)
-
\widehat{f}_{\varepsilon}(t)}
\\\nonumber
&=&
\frac{1}{nb^2}
\sum_{i=1}^n
 (\widehat{\varepsilon}_i(\widehat{\theta})
 - \widehat{\varepsilon}_i(\theta_o))
K_3^{(1)}
\left(\frac{\widehat{\varepsilon}_i(\theta_o)-t}{b}\right)
\\\nonumber
&& +
\frac{1}{2nb^3}\sum_{i=1}^n
(\widehat{\varepsilon}_i(\widehat{\theta})
-\widehat{\varepsilon}_i(\theta_o))^2
K_3^{(2)}
\left(\frac{\widehat{\varepsilon}_i(\theta_o) +
\beta(\widehat{\varepsilon}_i(\widehat{\theta})
-\widehat{\varepsilon}_i(\theta_o))-t}{b} \right),
 \label{rest1}
\end{eqnarray}
for some $\beta\in (0,1)$. In what follows, we will show that each
of the terms above is $o_{\prob}((nb)^{-1/2})$. First consider the
last term of (\ref{rest1}).  Since $\Lambda_{\theta}(y)$ and
$\widehat{m}_{\theta}(x)$ are both twice continuously
differentiable with respect to $\theta$, the second order Taylor
expansion gives, for some $\theta_1$ between
$\theta_o$ and $\widehat{\theta}$ (to simplify the notations, we assume
here that $p=\mbox{dim}(\theta)=1$),
\begin{eqnarray*}
\lefteqn{
 \widehat{\varepsilon}_i(\widehat{\theta})
 - \widehat{\varepsilon}_i(\theta_o)} \\
&=&
\Lambda_{\widehat{\theta}}(Y_i) - \Lambda_{\theta_o}(Y_i)
-\left(\widehat{m}_{\widehat{\theta}}(X_i) - \widehat{m}_{\theta_o}(X_i) \right)
\\
&=&
(\widehat{\theta}-\theta_o)
(\dot{\Lambda}_{\theta_o}(Y_i)  -  \dot{\widehat{m}}_{\theta_o}(X_i))
 +  \frac{1}{2} (\widehat{\theta}-\theta_o)^2 (\ddot{\Lambda}_{\theta_1}(Y_i)
  - \ddot{\widehat{m}}_{\theta_1}(X_i)).
\end{eqnarray*}
Therefore, since $\widehat{\theta}-\theta_o
=o_{\prob}((nb)^{-1/2})$
by Theorem 4.1 in Linton, Sperlich and Van Keilegom (2008)
(as before, we work with a slower rate than what is shown in the latter paper,
 since this leads to weaker conditions on the bandwidths), assumptions (A2) and (A7) imply that
\begin{eqnarray*}
\frac{1}{nb^3}\sum_{i=1}^n
(\widehat{\varepsilon}_i(\widehat{\theta})
-\widehat{\varepsilon}_i(\theta_o))^2
K_3^{(2)}
\left(\frac{\widehat{\varepsilon}_i(\theta_o)
+\beta(\widehat{\varepsilon}_i(\widehat{\theta})
-\widehat{\varepsilon}_i(\theta_o))-t}{b} \right)
=
 O_{\prob} \left((nb^3)^{-1} \right),
\end{eqnarray*}
which is $o_{\prob}((n b)^{-1/2})$, since
$(nb^5)^{-1} = O(1)$ under {\rm (A2)}.  For the first term
of (\ref{rest1}), the decomposition of
$\widehat{\varepsilon}_i(\widehat{\theta})
-\widehat{\varepsilon}_i(\theta_o)$ given above yields
\begin{eqnarray}
\nonumber
\lefteqn{
\frac{1}{nb^2}\sum_{i=1}^n
(\widehat{\varepsilon}_i(\widehat{\theta})
 -
 \widehat{\varepsilon}_i(\theta_o))
 K_3^{(1)}
 \left(\frac{\widehat{\varepsilon}_i(\theta_o)-t}{b} \right)
 }
\\\nonumber
&=&
\frac{(\widehat{\theta}-\theta_o)}
{nb^2}
\sum_{i=1}^n
(\dot{\Lambda}_{\theta_o}(Y_i)
-\dot{\widehat{m}}_{\theta_o}(X_i))
K_3^{(1)}
\left(\frac{\widehat{\varepsilon}_i(\theta_o)-t}{b}\right)
 +
o_{\prob}((n b)^{-1/2})
\\
&=&
\frac{(\widehat{\theta}-\theta_o)}
{nb^2}
\sum_{i=1}^n(\dot{\Lambda}_{\theta_o}(Y_i) -\dot{m}_{\theta_o}(X_i))
K_3^{(1)}\left(\frac{\varepsilon_i-t}{b}\right)
 +
 o_{\prob}((n b)^{-1/2}),
 \label{rest2}
\end{eqnarray}
where the last equality follows from a Taylor expansion applied to
$K_3^{(1)}$, the fact that
$$
\dot{\widehat{m}}_{\theta_o}(x) - \dot{m}_{\theta_o}(x)
=
 O_{\prob} ((nh)^{-1/2}(\log h^{-1})^{1/2}),
$$
uniformly in $x\in\mathcal{X}_0$ by Lemma \ref{mtheta}, and the
fact that $nhb^3 (\log h^{-1})^{-1} \rightarrow \infty$ under
(A2). Further, write
 \begin{eqnarray*}
\lefteqn{
\esp
\left[\sum_{i=1}^n
(\dot{\Lambda}_{\theta_o}(Y_i)-\dot{m}_{\theta_o}(X_i))
K_3^{(1)}\left(\frac{\varepsilon_i-t}{b}\right)
\right]}
\\
&=&
 \sum_{i=1}^n
 \esp
 \left[
 \dot{\Lambda}_{\theta_o}(Y_i)
  K_3^{(1)}
 \left(\frac{\varepsilon_i-t}{b}\right)
 \right]
 -
 \sum_{i=1}^n
 \esp
 \left[\dot{m}_{\theta_o}(X_i) \right]
 \esp
 \left[
 K_3^{(1)}\left( \frac{\varepsilon_i-t}{b} \right)
 \right]
 \\
&=&
A_n-B_n.
\end{eqnarray*}
We will only show that the first term above is $O(nb^2)$ for any
$t\in\Rit$. The proof for the other term is similar. Let
$\varphi(x,t)=\dot{\Lambda}_{\theta_o}(\Lambda_{\theta_o}^{-1}(m(x)+t))$ and set $\phi(x,t)
=\varphi(x,t)f_{\varepsilon}(t)$. Then, applying a Taylor expansion to
$\phi(x,\cdot)$, it follows that (for some $\beta \in (0,1)$)
\begin{eqnarray*}
\lefteqn{ A_n
= \sum_{i=1}^n
\esp
\left[
\dot{\Lambda}_{\theta_o}
 \left(\Lambda_{\theta_o}^{-1}
 (m(X_i)+\varepsilon_i)  \right)
 K_3^{(1)}
 \left( \frac{\varepsilon_i-t}{b}\right)
 \right]}
 \\
&=&
 n\int\int
 \phi(x,e)
 K_3^{(1)}\left(\frac{e-t}{b}\right)
 f_X(x) dx de
 \\
&=&
 n b
 \int\int
 \phi(x,t+bv)
 K_3^{(1)}(v) f_X(x)
  dx  dv
  \\
&=&
n b
\int \int
\left[
\phi(x,t)  +  b v \frac{\partial\phi}{\partial t}
(x,t+\beta b v) \right]
 K_3^{(1)}(v)
 f_X(x)
 dx  dv
 \\
 &=&
n b^2
\int  \int
v\frac{\partial\phi}{\partial t}(x,t+\beta b v)
K_3^{(1)}(v)
f_X(x)
dx  dv,
\end{eqnarray*}
since $\int K_3^{(1)}(v) dv=0$, and this is bounded by $Knb^2
\sup_{s:|t-s| \le \delta} \esp
|\frac{\partial
\phi}{\partial s} (X,s)| = O(nb^2)$ by assumption (A9)(ii). Hence,
Tchebychev's inequality ensures that
\begin{eqnarray*}
 \lefteqn{
 \frac{(\widehat{\theta}-\theta_o)}{b^2}
 \sum_{i=1}^n
 (\dot{\Lambda}_{\theta_o}(Y_i)
 - \dot{m}_{\theta_o}(X_i))
 K_3^{(1)} \left(\frac{\varepsilon_i-t}{b}\right)
  }
  \\
&=&
 \frac{(\widehat{\theta}-\theta_o)}{nb^2}
 O_{\prob} (nb^2+(n b)^{1/2})
=
o_{\prob}((nb)^{-1/2}),
\end{eqnarray*}
since $n b^{3/2} \rightarrow\infty$ by (A2). Substituting this in (\ref{rest2}), yields
$$
\frac{1}{nb^2} \sum_{i=1}^n
(\widehat{\varepsilon}_i(\widehat{\theta})
 - \widehat{\varepsilon}_i(\theta_o)) K_3^{(1)}
  \left( \frac{\widehat{\varepsilon}_i(\theta_o)-t}{b} \right)
 = o_{\prob}((nb)^{-1/2}),
$$
for any  $t\in\Rit$. This completes the proof. \hfill $\Box$ \\[-.3cm]

\noindent
{\bf Proof of Theorem \ref{Normality}.} It follows from Theorem \ref{Asymp-Dev} that
\begin{eqnarray}
\widehat{f}_{\widehat{\varepsilon}}(t)
- f_{\varepsilon}(t)
=
[\widetilde{f}_{\varepsilon}(t)
- \esp\widetilde{f}_{\varepsilon}(t)]
+
[ \esp
\widetilde{f}_{\varepsilon}(t) - f_{\varepsilon}(t)]
+
o_{\prob}((n b)^{-1/2}).
\label{asymexp1}
\end{eqnarray}
The first term on the right hand side of (\ref{asymexp1}) is
treated by Lyapounov's Central Limit Theorem (LCT) for
triangular arrays (see e.g. Billingsley 1968, Theorem 7.3). To
this end, let
$$
\widetilde{f}_{in}(t)
 =
 \frac{1}{b}
 K_{3}\left(\frac{\varepsilon_i-t}{b}\right).
$$
Then, under {\rm (A1), (A2)} and {\rm (A5)} it can be easily shown that
\begin{eqnarray*}
\frac{\sum_{i=1}^n
\esp\left|\widetilde{f}_{in}(t)-\esp
\widetilde{f}_{in}(t) \right|^3 }{\left(\sum_{i=1}^n \Var
\widetilde{f}_{in}(t) \right)^{3/2}}
 \leq
 \frac{C n b^{-2}f_{\varepsilon}(t)
 \displaystyle{\int}
\left| K_3(v) \right|^3 dv +
 o\left( nb^{-2} \right)}
 {
\left( n b^{-1}f_{\varepsilon}(t)
\displaystyle{\int}K_3^2(v) dv
+
 o\left(n b^{-1}\right)
\right)^{3/2}}
 =O((nb)^{-1/2})= o(1),
\end{eqnarray*}
for some $C>0$.  Hence, the LCT ensures that
$$
\frac{
\widetilde{f}_{\varepsilon}(t)
- \esp \widetilde{f}_{\varepsilon}(t)}
{\sqrt{\Var \widetilde{f}_{\varepsilon}(t)}}
=
 \frac{\widetilde{f}_{\varepsilon}(t)
 -\esp \widetilde{f}_{\varepsilon}(t)}
 {\sqrt{\frac{\Var\widetilde{f}_{1n}(t)}{n}}}
 \stackrel{d}{\rightarrow}
 N\left(0,1\right).
$$
This gives
\begin{eqnarray}
\sqrt{n b}
 \left(
 \widetilde{f}_{\varepsilon}(t)- \esp
\widetilde{f}_{\varepsilon}(t) \right)
\stackrel{d}{\rightarrow}
N\left(0, f_\varepsilon(t) \int K_3^2(v) dv \right).
\label{asymexp2}
\end{eqnarray}
For the second term of (\ref{asymexp1}), straightforward calculations show that
\begin{eqnarray*}
\esp \widetilde{f}_{\varepsilon}(t)- f_{\varepsilon}(t)
=
\frac{b^{q_3}}{q_3!} f_{\varepsilon}^{(q_3)}(t) \int v^{q_3}
K_3(v)  dv +  o(b^{q_3}).
\end{eqnarray*}
Combining this with (\ref{asymexp2}) and (\ref{asymexp1}), we
obtain the desired result. \hfill $\Box$ \\[-.3cm]

\begin{lem} \label{mtheta}
Assume {\rm (A1)-(A5)} and {\rm (A7)}. Then,
\begin{eqnarray*}
\sup_{x\in \mathcal{X}_0}
|\widehat{m}_{\theta_o}(x)
-m_{\theta_o}(x)|
 &=&
  O_{\prob}((nh)^{-1/2}(\log h^{-1})^{1/2}),
 \\
\sup_{x\in \mathcal{X}_0}
|\dot{\widehat{m}}_{\theta_o}(x)
- \dot{m}_{\theta_o}(x)|
 &=&
   O_{\prob}((nh)^{-1/2}(\log h^{-1})^{1/2}).
\end{eqnarray*}
\end{lem}

\vskip 0.3cm
 \noindent {\bf Proof.} We will only show the proof for
$\dot{\widehat{m}}_{\theta_o}(x) - \dot{m}_{\theta_o}(x)$, the
proof for $\widehat{m}_{\theta_o}(x) -m_{\theta_o}(x)$ being very
similar.
 Let $c_n=(nh)^{-1/2}(\log
h^{-1})^{1/2}$, and define
$$
\dot{\widehat{r}}_{\theta_o}(x)
 =  \frac{1}{nh}
 \sum_{j=1}^n  \dot{\Lambda}_{\theta_o}(Y_j)
 K_1\left( \frac{X_j-x}{h}  \right),
 \quad \dot{\overline{r}}_{\theta_o}(x)
  =  \esp  [\dot{\widehat{r}}_{\theta_o}(x)],
   \quad  \overline{f}_X(x)
 =  \esp[\widehat{f}_X(x)],
$$
where $\widehat{f}_X(x) =(nh)^{-1} \sum_{j=1}^n
K_1(\frac{X_j-x}{h})$. Then,
\begin{eqnarray}
\sup_{x\in \mathcal{X}_0}
|\dot{\widehat{m}}_{\theta_o}(x)  - \dot{m}_{\theta_o}(x)|
\leq
 \sup_{x\in \mathcal{X}_0}
 \left|
 \dot{\widehat{m}}_{\theta_o}(x)
 -  \frac{\dot{\overline{r}}_{\theta_o}(x)}
 {\overline{f}_X(x)} \right|
 + \sup_{x\in \mathcal{X}_0}
 \frac{1}{\overline{f}_X(x)}
 \left| \dot{\overline{r}}_{\theta_o}(x)
 -  \overline{f}_X(x)  \dot{m}_{\theta_o}(x)
  \right|.
 \label{mthetad}
\end{eqnarray}
 Since
$\esp[\dot{\Lambda}_{\theta_o}^4(Y)|X=x]<\infty$ uniformly in
$x\in\mathcal{X}$  by assumption (A7), a similar proof as was
given for Theorem 2 in Einmahl and Mason (2005) ensures that
$$
\sup_{x\in \mathcal{X}_0}
 \left| \dot{\widehat{m}}_{\theta_o}(x)
 -
 \frac{\dot{\overline{r}}_{\theta_o}(x)}
 {\overline{f}_X(x)}
 \right|
 =
 O_{\prob}
 \left(c_n\right).
 $$
 Consider now the second  term of (\ref{mthetad}).
 Since $\esp[\dot{\varepsilon}(\theta_o)|X]=0$, where
 $\dot{\varepsilon}(\theta_o)
 = \frac{d}{d\theta}(\Lambda_\theta(Y)-m_\theta(X))|_{\theta=\theta_o}$,
 we have
 \begin{eqnarray*}
\dot{\overline{r}}_{\theta_o}(x)
& = &
 h^{-1}
 \esp
 \left[
 \left\lbrace  \dot{m}_{\theta_o}(X)
 +  \dot \varepsilon(\theta_o)
 \right\rbrace
 K_1\left(  \frac{X-x}{h}
 \right) \right]
 \\
&=&
 h^{-1}
 \esp
 \left[
 \dot{m}_{\theta_o}(X)
 K_1\left(  \frac{X-x}{h}  \right)
 \right] \\
& = &
\int \dot{m}_{\theta_o}(x+hv)
K_1(v)f_X(x+hv) dv,
\end{eqnarray*}
from which it follows that
\begin{eqnarray*}
 \dot{\overline{r}}_{\theta_o}(x)
  -
  \overline{f}_X(x)
  \dot{m}_{\theta_o}(x)
= \int
\left[
\dot{m}_{\theta_o}(x+hv)
- \dot{m}_{\theta_o}(x)
\right]
K_1(v)
f_X(x+hv)
dv.
\end{eqnarray*}
Hence, a Taylor expansion applied to $\dot{m}_{\theta_o}(\cdot)$ yields
\begin{eqnarray*}
 \sup_{x\in\mathcal{X}_0}
  \left|
   \dot{\overline{r}}_{\theta_o}(x)
   -  \overline{f}_X(x)\dot{m}_{\theta_o}(x)
   \right|
=  O(h^{q_1})
 =
  O\left(c_n\right),
 \end{eqnarray*}
since $nh^{2q_1+1}(\log h^{-1})^{-1}=O(1)$ by (A2).   This proves that the second term
of (\ref{mthetad}) is $O(c_n)$, since it can be easily shown that $\overline{f}_X (x)$ is
 bounded away from $0$ and infinity, uniformly in
 $x\in\mathcal{X}_0$, using  (A3)(ii). \hfill $\Box$ \\[-.3cm]

\end{document}